\documentclass[runningheads]{cl2emult}

\usepackage{subeqnar} 
\usepackage{cropmark} 
\usepackage{math}     
\usepackage{cropmark} 

\usepackage{amsmath}   
\usepackage{amsfonts}
\usepackage{makeidx}

\newcommand{\al}{\alpha}

\newcommand{\C}{\mathbb{C}}
\newcommand{\Hy}{\mathbb{H}}

\newcommand{\T}{\mathbb{T}}

\newcommand{\Sumq}{\sum_{q\in\N}}

\newcommand{\diag}{\operatorname{diag}}
\newcommand{\gf}{geometrically finite}

\newcommand{\dpt}{\mathfrak{p}}

\newcommand{\lag}{\lambda_g}

\newcommand{\e}{\varepsilon}
\newcommand{\om}{\omega}
\newcommand{\ta}{\widetilde{a}}

\newcommand{\vphi}{\varphi}

\newcommand{\vwa}{very well approximable}

\newcommand{\cD}{{\cal{D}}}
\newcommand{\cE}{{\cal{E}}}

\newcommand{\jb}{\mathbf{j}}

\newcommand{\Be}{Besicovitch}
\newcommand{\BA}{\mathfrak{B}}

\newcommand{\Spr}{Sprind{\v z}uk}

\newcommand{\ip}{\,{\cdot}\,}

\newcommand{\ie}{{\it{i.e.}},\,}

\newcommand{\HD}{Hausdorff dimension}

\newcommand{\DA}{Diophantine approximation}
\newcommand{\N}{\mathbb{N}}
\newcommand{\Q}{\mathbb{Q}}
\newcommand{\R}{\mathbb{R}}
\newcommand{\Z}{\mathbb{Z}}
\newcommand{\ci}{\mathbb{S}^1} 
\newcommand{\col}{\colon}

\newcommand{\q}{\mathbf{q}}
\newcommand{\tv}{\mathbf{t}}

\newcommand{\dm}{dimension}

\newcommand{\Kh}{Khintchine}
 
\newcommand{\K}{{\cal{K}}}

\newcommand{\Ja}{Jarn{\'\i}k}
\newcommand{\JB}{Jarn{\'\i}k-Besicovitch}

\renewcommand{\d}{\delta}

\font\tenscr=rsfs10  scaled \magstep1
\font\sevenscr=rsfs7  scaled \magstep1
\font\fivescr=rsfs5  scaled \magstep1
\skewchar\tenscr='177 \skewchar\sevenscr='177 \skewchar\fivescr='177
\newfam\scrfam \textfont\scrfam=\tenscr \scriptfont\scrfam=\sevenscr
\scriptscriptfont\scrfam=\fivescr

\makeindex
\begin{document}
\title*{Exceptional sets in dynamical systems 
and Diophantine approximation} 
\titlerunning{Exceptional sets}
\author{M.~M.~Dodson}
\authorrunning{M. M. Dodson}

\institute{Department of Mathematics, 
  University of York, York YO10~5DD, UK} 

\maketitle

\begin{abstract}
\label{sec:ab}

The nature and origin of exceptional sets associated with the rotation
number of circle maps, Kolmogorov-Arnol'd-Moser theory on the
existence of invariant tori and the linearisation of complex
diffeomorphisms are explained.  The metrical properties of these
exceptional sets are closely related to fundamental results in the
metrical theory of Diophantine approximation.  The counterpart of
Diophantine approximation in hyperbolic space and a dynamical
interpretation which led to the very general notion of `shrinking
targets' are sketched and the recent use of flows in homogeneous
spaces of lattices in the proof of
the Baker-\Spr{} conjecture is described briefly.

\end{abstract}
\section{Introduction}
\label{sec:intro}

An important consequence of the natural connection between resonance
and Diophantine equations is the frequent occurrence of Diophantine
conditions in the study of dynamical systems.  A set in a space is
called {\em exceptional} if its members fail to satisfy some property
satisfied by `almost all' points in the space.  More precisely, let
$(X,\mu)$ be a measure space.  The set $E\subset X$ is exceptional if
$\mu(E)=0$.  From this point of view, the rationals $\Q$ form an
exceptional set in $\R$ and the Cantor set in $[0,1]$ is one with
respect to Lebesgue measure. The sets to be considered are general
`limsup' sets; the standard way of showing them to
be null is to use the easy half (convergence case) of the
Borel-Cantelli theorem.

As is well known, Hausdorff measure and dimension provide a means of
distinguishing sets of Lebesgue measure zero (or {\em null}\,)
sets.  In this article, the Hausdorff measure and dimension of some
exceptional sets of importance in dynamical systems will be discussed. 
The purpose of this paper is to show how these sets arise from 
certain Diophantine conditions failing to hold. The paper is not
intended to cover all recent developments which can be found in the
references cited.  We begin with a simple example; further details on
this and dynamical systems in general can be found in the
comprehensive introduction to dynamical systems of Katok and
Hasselblatt~\cite{KatokHass} or the more elementary text of Arrowsmith
and Place~\cite{APIDS}.

\section{Rotation number}
\label{sec:rotno} 
The rotation number\index{rotation number} is a measure of how far on
average a continuous orientation preserving homeomorphism $f\col
\ci\to \ci$ moves a point round the circle and illustrates nicely the
way Diophantine properties can arise in analysis
(see~\cite[Chapter~11]{KatokHass} or~\cite[Chapter~1]{NiteckiDD} for
accounts of this topic).  When $f$ is continuous it has a continuous
lift $F\col \R\to \R$ which satisfies
$ 
  e\circ F = f\circ e,
$ 
where $e\col\R\to \ci$ is the exponential map given by
\begin{equation*}
z=  e(t)=e^{2\pi i t}.
\end{equation*} 
The lift $F$ is unique up to translation by an integer. Moreover when
the homeomorphism $f$ is orientation preserving, then $f$ is of degree
1, $F$ is strictly increasing and the function $F-\text{id}$, given
by $F(t)-t$, is 1-periodic. Thus $F(t+1)-t-1=F(t)-t$ and
$F(t+1)=F(t)+1$.  When $f$ has a fixed point at $z_0$ say, so that $
f(z_0)=z_0$, then $F(t_0)=t_0+k_0$ for some $k_0\in\Z$, whence for each
$N\in \N$, the set of positive ($\ge 1$) integers, 
\begin{equation*}
F^N(t_0)=t_0+Nk_0,
\end{equation*}
and 
\begin{equation*}
\lim_{N\to\infty} \frac{F^N(t_0)}{N} = \lim_{N\to\infty}
\frac{t_0+Nk_0}{N}=k_0\in\Z,
\end{equation*}
where $F^N$ is the $N$-fold iteration $F\circ\dots \circ F$.  A
similar argument shows that the limit $\rho(F)$ for a $q$-periodic
point is a rational with denominator $q$. 
 The limit 
\begin{equation*}
\rho(F)=\lim_{N\to\infty} \frac{F^N(t)}{N}=\lim_{N\to\infty} \frac{F^N(t)-t}{N}
\end{equation*}
turns out to be independent of $t$ (for proofs
see~\cite[Prop.~11.1.1]{KatokHass} and~\cite[p~33]{NiteckiDD}).  In
fact it can be shown that $\rho(F)$ is irrational if and only if $f$
has no periodic points~\cite{KatokHass,NiteckiDD}.

  Two rotation numbers $\rho(F)$ and $\rho(F')$ for lifts $F,F'$ for $f$
differ by an integer and so the rotation number $\rho(f)$ for $f$ can
defined 
as the fractional part $\{\rho(F)\}$ of $\rho(F)$. As a simple example, 
the rotation number of a rotation $r_\alpha$
by an angle $2\pi\alpha$, where $0\leq \alpha <1$,  given by
  \begin{equation}
\label{eq:rotfn}
    r_\alpha(z)=z\,e(\alpha) = ze^{2\pi i \alpha}  
  \end{equation}
  is, naturally enough,  $\alpha$. 
It follows immmediately that $\rho(f)$ is irrational if and only if
$f$ has no periodic points.  If $\rho(f)$ is irrational then for $z\in
\ci$, the closure $A$ of the orbit
$$
\om(z)=\{f^n(z)\col n\in\N\}
$$
does not depend on $z$ and $A$ is either perfect and nowhere dense
or $A=\ci$.  In the latter case $f$ is transitive and is topologically
conjugate\index{conjugacy} to the rotation $r_{\rho(f)}\,$, \ie there
exists an orientation preserving homeomorphism $\vphi\col \ci\to \ci$
such that
\begin{equation}
\label{eq:conj}
  f=\vphi^{-1}\circ r_{\rho(f)} \circ \vphi \quad {\text{ or }} 
\quad \vphi \circ f= r_{\rho(f)} \circ \vphi, 
\end{equation}
written $f\sim r_{\rho(f)}$.

Denjoy showed that when $f$ is $C^2$
rather than just $C^1$, the irrationality of $\rho(f)$ implies that
$f$ is transitive so that $f\sim r_{\rho(f)}$, \ie $f$ is
topologically conjugate to a rotation.  This partial converse to the
rotation $r_\rho$ having rotation number $\rho$ is best possible as
there are counterexamples when $f$ is $C^{2-\e}$~\cite{Yoccoz92}.

Two kinds of \DA{} enter the picture when the differentiability of the
conjugacy function $\vphi$ is considered.  When $\vphi$ is $C^r$,
where $r$ is finite or infinite (but $\vphi$ is not analytic), the
existence or otherwise of a smooth conjugacy is determined by the
Diophantine type of $\rho$ (see \S\ref{sec:DT} below). If $f$ is $C^r$
and $\rho(f)$ is of Diophantine type then $f$ is smoothly conjugate to
$r_{\rho(f)}$.  The precise details of the differentiability
conditions are rather complicated; further details and references are
in~\cite{KatokHass} and~\cite{Yoccoz92,Yoccoz95a}.

When  $\vphi$ is analytic and $\rho$ is a Bruno 
number\index{Bruno number}, 
\ie $\rho\in\R\setminus\Q$ and  the denominator $q_n$ of the 
$n$-th convergent in the continued fraction expansion of
$\rho\in \R\setminus \Q$ satisfies
\begin{equation*}
   \sum_{n=1}^\infty \frac{\log q_{n+1}}{q_n}<\infty,
\end{equation*}
then any analytic $f$ with rotation number $\rho$ is analytically
conjugate to $r_\rho$ (continued fractions and convergents are
explained in~\cite{Cassels,HW}).  This condition is also optimal since
if $\rho$ is not a Bruno number, there is a Blaschke product with
rotation number $\rho$ not analytically conjugate to
$r_\rho$~\cite{Yoccoz92}.

To illustrate the ideas, we shall consider analytic conjugation when
$f$ is close to a rotation and see how Diophantine conditions emerge.
Suppose the diffeomorphism $f$ is analytic with (analytic) lift $F$
given by
\begin{equation*}
F(x)=x+\rho + a(x),
\end{equation*}
where $a$ is a 1-periodic analytic function, real on the real axis
and bounded on  the strip $\{z\in\C\col |\Im z|< \e\}$ for some $\e>0$. 
 Consider the lifts $F$, $\Phi$ and $R_\rho$ of 
$f$, $\vphi$ and $r_\rho$ respectively and the lift 
\begin{equation}
\label{eq:lift2}
  F \circ \Phi(t) = \Phi \circ R_\rho(t) 
\end{equation}
of equation~\eqref{eq:conj}, where $R_\rho$ is translation by $\rho$
given by $R_\rho(t)=t+\rho$ and where $\Phi$ is assumed to be close to
the identity, \ie $\Phi(t)=t+h(t)$, where $h(t)$ is small.
Then~\eqref{eq:lift2} reduces to the functional equation
\begin{equation*}
  t+ h(t)+ \rho +a(t+h(t)) = t+ \rho + h(t+\rho),
\end{equation*}
\ie to 
\begin{equation}
\label{eq:fnaleq}
 h(t+\rho) - h(t) = a(t+h(t)). 
\end{equation}
One then seeks a solution to this equation  by successive
approximation. 
Since $a$ and $h$ are small, consider the first order approximation
\begin{equation}
\label{eq:fofe}
   h^0(t+\rho)- h^0(t)=\ta(t)=a(t)-a_0,
\end{equation}
where $a_0=\int_0^1 a(t)dt$, the mean of $a(t)$.  
Since $a$ and $h$ are periodic, they have Fourier series expansions 
\begin{equation*}
  a(t)=\sum_{k\in\Z } a_ke^{2\pi i kt} \quad \text{ and } 
  h^0(t)=\sum_{k\in\Z } h^0_ke^{2\pi i kt}.
\end{equation*}
On substituting in~\eqref{eq:fofe} and comparing the coefficients of
$e^{2\pi i kt}$, the coefficients of the linear approximation $h^0$ to
the unknown function $h$ can be determined in terms of those of the
function $a$, as follows: $h^0_0=0$ and for $k\ne 0$,
\begin{equation}
\label{eq:h0k}
  h^0_k=\frac{a_k}{e^{2\pi i k\rho}-1}.
\end{equation}
It is evident that for $h^0$ to exist, the denominator and numerator
of $h^0_k$ must vanish together (this is why the function
$\ta(t)=a(t)-a_0$ is introduced).  Since $a$ is analytic, its Fourier
coefficients $a_k$ decay rapidly and 
$$
|a_k|<\d\,e^{-|k|\e},
$$ 
where $\sup\{|a(z)|\col |\Im z|<\e\}<\d $.  
Now as $2\theta/\pi\leq \sin \theta \leq \theta$ when $0\leq
\theta\leq \pi/2$, it is readily verified that for any real $\theta$,
$$ 
|e^{i\theta}-1|=2 |\sin \theta/2| \geq 2 \| \theta\|/\pi,
$$
where $\| \theta\| = \theta$ when $0\leq \{\theta\} \leq
1/2$ and $1-\{\theta\}$ otherwise ($ \{\theta\}$ is the fractional
part of the real number $\theta$). 
Hence when $k\ne 0$,
\begin{equation*}
|e^{2\pi i k\langle\rho \rangle} -1|= |e^{2\pi i k\rho}-1| 
\geq 2\| k\rho\|/\pi.
\end{equation*}
Note that by~\eqref{eq:h0k}, the denominator of $h^0_k$ is $e^{2\pi i k\rho}-1$.

Dirichlet's theorem~\cite[\S11.2]{HW}\index{Dirichlet's theorem}
states that given any real number $\rho$ and positive integer $N$,
there exists a positive integer $k\leq N$ and an integer $j$ such that
\begin{equation*}
  \left|\rho-\frac{j}{k}\right|< \frac1{kN}\leq \frac1{k^2}.
\end{equation*}
It follows that there exist infinitely many $k\in\N$ such that 
\begin{equation*}
  \|k\rho\| \leq k^{-1}
\end{equation*}
(recall that $ \|k\rho\|=\min\{|k\rho-p|\col p\in\Z\}$).  Indeed for
rational $\rho$, the denominator will vanish for infinitely many $k$,
so destroying the convergence of the Fourier series for $h$. Even for
irrational $\rho$, the denominator will get arbitrarily small and so
make convergence problematic.  However, the very rapid decay in the
numerator $a_k$ can be set against this to obtain convergence for
those numbers $\rho$ for which $\|k\rho\|$ does not get too small.  To
be more precise, comparing the moduli of the denominator and numerator
of $h^0_k$ in~\eqref{eq:h0k}, we see that the Fourier series for $h^0$
will converge if $\rho$ satisfies
\begin{equation}
  \label{eq:DTcondn}
 \|k \rho\|  \geq K|k|^{-v}
\end{equation}
for some $K,v>0$.  By Dirichlet's theorem, we must take $v>1$.

\subsection{Diophantine type}\index{Diophantine type}
\label{sec:DT}

A number $\rho$ satisfying 
\begin{equation*}
  \left|\rho-\frac{j}{k}\right|\geq \frac{K}{k^{\sigma+2}},
\end{equation*}
for some $\sigma>0$ and for all $k\in \N$, is said to be of {\em Diophantine
  type} $(K,\sigma)$~\cite{ArnoldMM} and the set of such numbers is
denoted by $\cD(K,\sigma)$; the set of numbers of Diophantine type
$(K,\sigma)$ for some positive $K,\sigma$ is denoted
$\cD=\bigcup_{K,\sigma} \cD(K,\sigma)$.  Similar conditions turn up in
other settings, as we will see below.  The set of Bruno numbers
includes $\cD$ since if $\xi\in\cD(K,\sigma)$, then
$q_{n+1}=O(q_n^{\sigma})$ and $ \sum_n (\log q_{n+1})/q_n <\infty$.
Note that ${\cD}(0)=\bigcup_{K>0} \cD(K,0) =\BA$, the set of badly
approximable numbers\index{badly approximable}.

Using a version of Newton's tangent method, a sequence of successive
approximations $\Phi_n$ is constructed which converges to an analytic
function $\Phi$ satisfying~\eqref{eq:lift2}, again providing $\rho$ is
of type $(K,\sigma)$ for some $K$, $\sigma$ and $a$ is sufficiently
small. Projecting from the lift, the
desired conjugacy exists, {\em providing} $\rho$ is of type
$(K,\sigma)$, \ie providing $\rho\in \cD$.
More details are in~\cite[\S12]{ArnoldGM}.

The size of the complementary set $E=\R\setminus\cD$ of points for
which the iterative method does not guarantee convergence -- and so
conjugacy -- should be as small as possible.  The exponent $\sigma$ in
the Diophantine type affects the range of validity of the result
significantly, so we consider the set $E_\sigma$ given by
\begin{eqnarray*}
    E_\sigma = \bigcup_{K>0}\{\xi\in{\R} 
\col |\xi-p/q|< Kq^{-2-\sigma} {\text{ for some }}  p/q &\in{\Q}\}\\ 
             = \bigcup_{K>0} \{\xi\in{\R} 
\col \|q\xi\|< Kq^{-1-\sigma} {\text{ for some }} q\in{\N}\}.&  
\end{eqnarray*} 
The set 
\begin{equation*}
  E=\bigcap_{\sigma>0} E_\sigma =\bigcap_{\sigma\in\Q^+} E_\sigma 
=\lim_{\sigma\to \infty} E_\sigma
\end{equation*}
consists  of points not of Diophantine type for any $K,\sigma>0$.  
 Given $\sigma>0$,  almost all real numbers are
 of Diophantine type $(K,\sigma)$ for some positive $K$, 
\ie  the set 
 \begin{equation}
\label{eq:Dsigma}
 {\cD}(\sigma) =  
\bigcup_{K>0} \{\xi\in{\R}\col |\xi-p/q|\geq Kq^{-2-\sigma} {\text{
       for each }} p/q\in{\Q}\} 
 \end{equation}
is of full measure,  since the complementary set 
 \begin{equation*}
  E_{\sigma}= \bigcap_{K>0} \{\xi\in\R\col |\xi-p/q|< Kq^{-2-\sigma}
  {\text{ for   some }}  p/q\in\Q\} 
 \end{equation*}
 is null.  This can be proved directly but we will consider a related
 set $\K_v$ and use Khintchine's theorem.  This theorem also implies
 that the set $\BA$ of badly approximable numbers is null. Note that
 badly approximable numbers are good for conjugacy.

\subsection{Very well approximable numbers 
and Khintchine's theorem}\index{Khintchine's theorem}

The set of $v$-{\em approximable }\index{$v$-approximable} numbers
is defined as
 \begin{equation}
   \label{eq:Kv}
 \K_v =
\{\xi\in\R\col \|q\xi\|<q^{-v} {\text{ for infinitely many }} q\in\N \}; 
 \end{equation}
 when $v>1$, the points are often called very well approximable
 (VWA)\index{VWA (very well approximable)}.  The set $ \K_v$ can be
 usefully expressed as a `limsup' set\index{limsup set}, \ie
\begin{equation} 
  \label{eq:limsup}
   \K_v=\bigcap_{N=1}^\infty \bigcup_{q=N}^\infty B_{\psi(q)}(q) 
= \limsup_{N\to\infty}  B_{\psi(N)}(N), 
\end{equation}
where $B_\delta(q)=\{\xi\in\R\col \|q\xi\|<\delta\}$~\cite{FalcGFS}.
 It is not difficult to verify that for any $\e>0$,
\begin{equation}
  \label{eq:KEK}
  \K_{\sigma + 1+\e}\subset E_\sigma \subset \K_{\sigma+1},
\end{equation}
see~\cite[\S7.2]{MDAM}.  The set $\cD(\sigma)$ (defined
in~\eqref{eq:Dsigma}) is complementary to the sets of points
approximable to exponent $\sigma+2$ and consists of points which are,
in a sense, not close to the rationals or \vwa{} points.  Points in
$\cD(\sigma)$ have desirable approximation properties for certain
applications.

Khintchine's theorem relates the Lebesgue measure of sets more general
than $\K_v$ to the convergence or divergence of a sum.  The Lebesgue
measure of a subset $A$ of Euclidean space will be denoted by $|A|$.
A set $A$ is {\em full} if the Lebesgue measure of its complement is null.
First let $\psi:\N \to \R^+$ be a function and let
\begin{equation}
  \label{eq:Kpsi}
\K(\psi)=\{\xi\in\R\col\|q\xi\|<\psi(q){\text{ for infinitely many }} 
q\in\N \} 
\end{equation}
 be the set of $\psi$-{\em approximable} points\index{$\psi$-approximable}. 
\begin{theorem}
 If the sum 
$ 
  \Sumq \psi(q)
$ 
converges,  then $\K(\psi)$ is null, while if the sum diverges
and $\psi$ is decreasing,  $\K(\psi)$  is full. 
\end{theorem}
\Kh's therorem is reminsicent of the Borel-Cantelli lemma and indeed
the proof when the sum converges is the same, as it
suffices to show that
\begin{equation*}
  \sum_q |\{\xi\in\R\col \|q\xi\|<\psi(q)\}|<\infty.
\end{equation*}
The divergent case is, however, much more difficult and relies on
pairwise (quasi-)independence and deeper ergodic ideas. Proofs of
Khintchine's theorem and generalisations are
in~\cite[Chapter~VII]{Cassels},~\cite{HarmanMNT} and~\cite{Sprindzuk};
the last two references include discussions of `zero-one' laws.

Since $\K(\psi)=\K_v $ when $\psi(x)=x^{-v}$ and since the sum $\sum_q
q^{-v}$ converges for $v>1$, Khintchine's theorem implies that $
|\K_{\sigma+1}|=0$ for $\sigma>0$.  It follows from~\eqref{eq:KEK}
that $|E_\sigma|\le |\K_{\sigma+1}|=0$ for $\sigma>0$.  Since
$\cD:=\bigcup_{\sigma>0}\cD(\sigma)= \lim_{\sigma\to
  \infty}\cD(\sigma)$, $\cD$ is of full measure and its complement,
the set of of numbers {\em not} of Diophantine type $\cD(K,\sigma)$
for any $K>0$, is an exceptional set.

There is a very general higher dimensional form of the theorem due to
Groshev which is useful to us.  Let $\K^{(m,n)}(\psi)$ be the set of
$m\times n$ real matrices $X$ such that
\begin{equation}
  \label{eq:KGineq}
\|\q \, X\|:= \max\{|\q\ip X_1\|,\dots,\|\q\ip X_n\|\} < \psi(|\q|), 
\end{equation}
where $X_i$ is the $i$-th column of $X$ and $|\q|=
\max\{|q_1|,\dots,|q_m|\}$ is the height of $\q$, for infinitely many
$\q\in\Z^m$.  The set of $m\times n$ real matrices is identified in
the natural way with $\R^{mn}$.
\begin{theorem}
 If the sum 
$ 
  \sum_{{\q}\ne 0} |\q|^{m-1}\psi(|\q|)^n
$ 
converges,  then $\K^{(m,n)}(\psi)$ is null, while if the sum diverges
and $\psi$ is decreasing,  $\K^{(m,n)}(\psi)$  is full. 
\end{theorem}

\subsection{Hausdorff dimension\index{Hausdorff dimension}  
and the \JB{} theorem\index{Jarn{\'\i}k-Besicovitch theorem} }  

The dependence of the exceptional sets $E_\sigma$ and $\K_{\sigma +
  1}$ on $\sigma$, where $\sigma>0$, is revealed through their
Hausdorff dimension (definitions and expositions are 
in~\cite{MDAM,FalcGFS,MattilaGS,Rogers}).  In fact the dimension
$\dim \K_{\sigma+1}$ of the latter is given by another famous result
in number theory, namely the \JB{} theorem~\cite{Bes34,Ja29}
(see~\cite{MDAM} for more details).
\begin{theorem}
  When $v>1$,
  \begin{equation*}
    \dim \K_v=\frac{2}{v+1}.
  \end{equation*}
\end{theorem}
By Dirichlet's theorem, when $v \le 1$, $\K_{v}=\R$, whence $ \dim
\K_{v}=1$~\cite{FalcGFS,MattilaGS}.  The proof that $2/(v+1)$ is an upper
bound for the \dm{} follows straightforwardly from $\K_v$ being a
limsup set and a natural covering argument. To prove it a lower
bound is much harder. \Ja's original proof and his more general
Hausdorff measure result for simultaneous Diophantine approximation
relied heavily on arithmetic arguments.  \Be's later independent proof
was more geometric and was a basis for the widely used regular
systems~\cite{BS} and ubiquitous systems~\cite{DRV90a}. 

It follows from the inclusions~\eqref{eq:KEK} and the \JB{} theorem
that when $\sigma\geq 0$,
\begin{equation*}
  \frac{2}{\sigma+2+\e} \leq \dim E_\sigma \leq
  \frac{2}{\sigma+2}.
\end{equation*}
Since $\e>0$ is arbitrary,  $\dim
E_\sigma=2/(\sigma+2)$.  Now let $E=\bigcap_{\sigma>0}
E_\sigma=\lim_{\sigma\to\infty} E_\sigma$ since $E_\sigma$ decreases
as $\sigma$ increases. Then
\begin{equation}\
\label{eq:dimE0}
\dim E=  \lim_{\sigma\to\infty} \frac{2}{\sigma+2} 
= 0,
\end{equation}
and the complement of $\cD$ has \HD{} 0.  Thus when $f$ is analytic,
$f$ is analytically conjugate to a rotation unless the rotation number
of $f$ lies in a set of \HD{} 0.  Rational points lie in this set and
they play a special role in this problem.  Because of the connection
with the physical phenomenon of resonance, the rationals are called
{\em resonant}. The set $E_\sigma$ can be thought of as consisting of
families of regularly spaced points and associated `fractal dust'.
More details are given in~\cite{ArnoldGM,NiteckiDD,Yoccoz92} and there
is a simplified account in~\cite{FalcGFS}.  The phenomenon is
widespread, see for example \S4 below on linearisation and in normal
forms~\cite{dgy95,DRV94}.

In higher dimensions, the \HD{} of the set 
\begin{equation}
  \label{eq:Knv}
\K_v^{(1,n)}=\{\xi\in\R^n\col\|q \xi\|
<q^{-v}{\text{ for infinitely many }} q\in\N \}  
\end{equation}
is $(n+1)/(v+1)$ when $v>1/n$~\cite{Ja31}; and that of the set 
\begin{equation}
  \label{eq:Gnv}
\K_v^{(n,1)}=\{\xi\in\R^n\col\|\q\ip \xi\|<|\q|^{-v}{\text{ for infinitely many }} \q\in\Z^n \}  
\end{equation}
is $n-1+(n+1)/(v+1)$ when $v>n$~\cite{BD86}. For the general systems
of linear forms which combines both results, see~\cite{mmd92}.
Levesley has proved a general inhomogeneous form~\cite{JL98} which was
used by Dickinson in normal forms for pseudo-elliptic
operators~\cite{dgy95}.

\Ja{} also showed that the exceptional set of badly approximable
numbers has full Hausdorff dimension 1~\cite{Ja28}.  This interesting
result was extended greatly by Schmidt~\cite{Schmidt69} who showed
that the set was `thick' and by Dani~\cite{Dani85,Dani86a} who
extended the ideas to dynamical systems. Kleinbock obtained a general
inhomogeneous version of \Ja's theorem by exploiting these
ideas~\cite{Kleinbock01,Kleinbock99}, for an expository account of
dimension and dynamical systems, see~\cite{Pesin97}.

\section{Kolmogorov-Arnol'd-Moser (KAM) theory\index{KAM theory}} 
\label{sec:KAM}
The next example comes from the study of the stability of the solar
system~\cite{Moser78}.  This is one of the oldest problems in
mechanics and is of course a special case of understanding the motions
of $N$ bodies of point mass subject to Newtonian attraction.  The
solution is well known for $N=2$ and periodic solutions exist, with
the bodies moving in elliptic orbits about their centre of mass.  In
the absence of any effects such as friction, the solution persists for
all time.  For $N\geq 3$, however, the situation is extraordinarily
complicated and far from understood.  This is the case for solar
systems, in which one of the bodies is a sun with mass $m_N$ very much
larger than the other masses $m_1,\dots,m_{N-1}$ of the planets.  If
as a first approximation, the centre of mass of the system is assumed
to coincide with that of the sun and if the gravitational interactions
between the planets and other effects are neglected, then the system
decouples into $N-1$ two-body problems, in which each planet describes
an elliptical orbit around the sun, with period $T_j$ say and
frequency $\omega_j=2\pi/T_j$, $j=1,\dots,N-1$.

For each vector $\omega=(\omega_1,\dots,\omega_n)$ of frequencies in
 the $n$-\dm al
 torus $\T^n=\ci\times \dots\times \ci$, 
the map $\varphi_{\omega} \col \R \to \T^n$  given by
\begin{equation}
\label{eq:qpflow}
  \varphi_{\omega}(t) = \varphi_{\omega}(0)+ t\omega
\end{equation}
is a quasi-periodic flow on the torus\index{quasi-periodic flow}.  The
case $n=1$ corresponds to uniform motion around a circle and so is
periodic.  When the frequencies are all rational, the flow is a closed
path on the torus and returns to its starting point after a time equal
to the lowest common multiple of the denominators of the frequencies.
If the frequencies are not all rational, then by Kronecker's
theorem~\cite{HW}\index{Kronecker's theorem}, the flow winds round the
torus, densely filling a subspace of dimension given by the number of
rationally independent frequencies. Thus when the frequencies are
independent over the rationals, the closure of $\varphi_{\omega}(\R)$
is the torus $\T^n$.  Solutions given by quasi-periodic flow
$\varphi_{\omega}$ on the torus $\T^n$, where $n=N-1$ will persist for
all time.

The planetary interactions are represented by a small perturbation of
the ori\-ginal Hamiltonian.  The stability of the solar system then
reduces to the question of whether the solutions of the perturbed
Hamiltonian system will continue to wind round a perturbed invariant
torus, so that motion of the planets remains quasi-periodic and
persists for all time.  This model is of course idealised and takes no
account of the long term fate of the universe.  Weierstrass
constructed a formal series solution but was unable to establish
convergence because the denominators in the terms of the series
contained factors of the form $\q\ip \omega$ which can become very
small for certain integer vectors $\q$~\cite[Chapter 1,
\S2]{MoserSRM}.  Investigating the three body problem led Poincar\'e
to speculate that quasi-periodic solutions need not exist. However,
Siegel overcame a related `small denominator' problem arising in the
linearisation of complex diffeomorphisms
(\cite{MoserSRM,Siegel42,Siegel52}, see~\S4 below) and a little later
Kolmogorov formulated the theorem that quasi-periodic solutions for a
perturbed analytic Hamiltonian system not only existed but were
relatively abundant in the sense that they formed a complicated Cantor
type set of positive Lebesgue measure~\cite{Kolmogorov54}.  This
theorem was proved completely in 1962 by Arnol'd~\cite{Arnold63};
independently Moser proved an analogous result for sufficiently
smooth, area-preserving planar maps (`twist'
maps)~\cite{Moser62,SiegelMoserCM},
the degree of differentiablility and the Diophantine conditions being
relaxed substantially by R\"ussmann~\cite{Russmann83,Russmann91} (see
also~\cite[\S6.3]{APIDS},\cite{NTDS},~\cite[Chap.~1]{MoserSRM},~\cite{Poschel82}).
It follows that for planets with mass very much less than the sun and
for the majority of initial conditions in which the orbits are close
to co-planar circles, distances between the bodies will remain
perpetually bounded, \ie the planets will never collide, escape or
fall into the sun.  Further details can be found
in~\cite{AKNIII,KAMtut}.

As in~\S2, the Diophantine condition~\eqref{eq:kamda} in the theorem
below guarantees convergence of the Fourier series and of an infinite
\dm al extension of Newton's tangent method.  The theorem is expressed
in terms of general analytic Hamiltonians and has been taken
from~\cite{MoserSRM}, p.~44.

\begin{theorem}[Kolmogorov-Arnol'd]
  \ Let $U$ be an open bounded subset of $ \R^n$ and let the function
  $H(x,y,\mu)$  be a real analytic function of $x,y,\mu$ for all $x,y\in U$ 
and near $\mu=0$. Moreover, $H$ is assumed to have period $2\pi$ in 
$x_1,\dots,x_n$ and $H_0=H(x,y,0)$ to be independent of $x$. 
Let $c\in\R^n$ be chosen so that 
  the frequencies 
$$
\omega_k=\frac{\partial H_0(c)}{\partial x_k}, \ \  1\leq k  \leq n, 
$$ 
in $\omega=(\omega_1,\dots,\omega_n)$ satisfy  the
Diophantine condition 
\vspace{.1in}
\begin{equation}
\label{eq:kamda}
|\q\ip \omega|=| q_1 \omega_1 + \dots + q_n \omega_n| \geq K|\q|_1^{-v}, 
\end{equation}
where $|\q|_1=|q_1|+\dots+|q_n|$, for some positive constants
$K=K(\omega)$ and $v=v(\omega)$  for all non-zero
$\q=(q_1,\dots, q_n)\in\Z^n$.  Suppose further that the Hessian 
$$
\det
\left(\frac{\partial^2 H_0(c)}{\partial x_j \partial x_k} \right)\ne 0.
$$ 
 Then for sufficiently small $|\mu|$, there exists
  an invariant torus described by
\begin{equation*}
  (x,y) = \left(\theta+ f(\theta,\mu), c + g(\theta,\mu)\right),
\end{equation*}
where  $\theta=(\theta_1,\dots,\theta_n)$ have period $2\pi$ in each
component and $f(\theta,\mu)$, $g(\theta,\mu)$ are analytic functions of $\mu$
which vanish for $\mu=0$.
Moreover the flow on this torus is given by
\begin{equation*}
\dot  \theta=\omega.
\end{equation*}
\end{theorem}
 The exponent $v$ is
 subject  to two conflicting requirements.  It
should be large enough ($v>n$) to ensure that the Diophantine 
condition above 
is not too restrictive but small enough to ensure that the perturbation
has physical significance and that the stability is robust.   
 The
proof breaks down when the frequencies lie in the complementary
exceptional set $E_v$ say of frequencies which are close to resonance
in the sense that given any $C>0$, there exists a $\q\in\Z^n$ such
that
\begin{equation*}
|\q\ip x|<C|\q|_1^{-v}. 
\end{equation*}
This set is  closely related to the set 
\newcommand{\hcL}{{\widehat{\cal{L}}}} 
\begin{equation}
\hcL_v(\R^n) =
\{x\in\R^n\col |\q\ip x|<|\q|^{-v} \text{ for infinitely many } \q\in\Z^n\}
  \label{eq:cL7}
\end{equation}
and in fact by an argument similar to that giving~\eqref{eq:KEK}, 
for any $\e>0$, 
\begin{equation*}
  \hcL_{v+\e}(\R^n)\subset E_v\subset \hcL_v(\R^n) 
\end{equation*}
(see~\cite[\S7.5.2]{MDAM}).  This inclusion implies that the two sets $
\hcL_v(\R^n)$ and $E_v$ have the same \HD.  The set $\hcL_v(\R^n) $ is
related to $\K^{(n,1)}_v$ and the \HD{} of $\hcL_v(\R^n) $ is a special case
of an analogue, proved by Dickinson~\cite{HD93}, of the general
form of the \JB{} theorem: 
\begin{equation}
\label{eq:dimhL}
\dim \hcL_v(\R^n)\,  (= \dim E_v)\, = n-1+ \frac{n}{v+1}.  
\end{equation}
when $v>n-1$ (see
also~\cite{DV}); note that $ \hcL_v(\R^n)=\R^n$ otherwise.

The resonance sets $R_{\q}=\{x\in\R^n\col \q\ip x=0\}$ have \dm{}
$n-1$ and $ \hcL_v(\R^n)$ consists of families of hyperplanes through
the origin with integer vector normals plus `fractal dust'.
\newcommand{\tcL}{\widetilde{\cal{L}}}
Another approach is via `averaging'~\cite{ArnoldGM}; this
can involve Diophantine approximation on manifolds~\cite{DRV89a}.

\section{Linearisation\index{linearisation}}
\label{sec:lin}

The linearisation of a complex analytic diffeomorphism $f:\C^{\,n}\to
\C^{\,n}$ in a neighbourhood of a fixed point is related to the
rotation number. The fixed point can be taken to be the origin and the
linearisation means that near the origin $f$ is analytically conjugate
to its linear part (Jacobian) $Df|_0=A$ say.  The linearising
transformation $\phi$ is given by the solution to the functional
equation
\begin{equation*}
f\circ \phi=\phi\circ A  
\end{equation*}
and when $n=1$ is known as Schr{\"o}der's
equation\index{Schr{\"o}der's equation}.  Linearisation can be 
regarded as obtaining a normal form and is analogous to diagonalising
a matrix.  Problems similar to those discussed above arise when the
eigenvalues $\al_1,\dots,\al_n$ of $A$ are close to being resonant in
the sense that they are close to satisfying the equation
\begin{equation*}
  \al_k=\prod_{r=1}^n \al_r^{j_r}
\end{equation*}
for some $k$ and all $\jb=(j_1,\dots,j_n)$ with $|\jb|_1\geq 2$ and
$j_r\geq 0$, $r=1,\dots,n$.  Linearisation is well understood when
$n=1$ and the diffeomorphism $f:\C\to \C$ with $f(0)=0$ can be
linearised when $|(Df)|_0|=|f'(0)| \not=1$.  The interesting case when
$|f'(0) |=1$ is closely related (via lifts) to the conjugacy of a
circle map to a rotation, discussed in~\S\ref{sec:rotno} above, and
necessary and sufficient conditions for the linearisation of a
diffeomorphism $f\col\C\to\C$ in terms of Bruno numbers are
known~\cite{Herman87,Yoccoz92,Yoccoz95a}.  Less is known when $n \geq 2$ 
 but Siegel established sufficient conditions on $Df|_0$ which
guarantee the existence of a linearising transformation $\phi$.  By 
analogy with the terminology Diophantine type, the point
$(\al_1,\dots,\al_n)$ in $\C^n$ is said to be of {\em multiplicative
  type} $(K,v)$~\cite[p.~191]{ArnoldGM} if
\begin{equation}
\left|\al_k-\prod_{r=1}^n \al_r^{j_r}\right|\geq K|\jb|_1^{-v}
  \label{eq:multype}\end{equation}
for all $\jb\in (\N\cup \{0\})^n$ with $|\jb|_1\geq 2$ and all
$k=1,\dots,n$.  Siegel showed that if the
vector $(\al_1,\dots,\al_n)$ of eigenvalues of $Df|_0$ is of
multiplicative type $(K,v)$ for some $K,v>0$,
then $f$ can be linearised locally~\cite{Siegel42,Siegel52}. 
Siegel's condition is not seriously restrictive when $v>(n-1)/2$ as
then for any $K>0$ the set of points of multiplicative type $(K,v)$
has full measure.  However, the size of the neighbourhood of
linearisation depends on $v$ (it also depends on $K$ but less
significantly) and so it is desirable not to make $v$ too large.

Given an exponent $v$, the complementary set of points ${\cE}_v$ say
in $\C^{\,n}$ (regarded as $\R^{2n}$) consists of points which fail to
be of multiplicative type $(K,v)$ for any $K>0$ and so fail to satisfy
the conditions of Siegel's theorem. 
The multiplicative Diophantine condition~\eqref{eq:multype} is hard to
handle and the question is simplified by means of the map
$$
(z_1,\dots,z_n)\mapsto (e^{2\pi iz_1},\dots,e^{2\pi iz_n})
$$ 
which preserves the \HD{} and reduces the problem 
to sets with a simpler structure~\cite{DRV94}. 
Multiplicative type is replaced by the  more amenable notion of a 
point $z=x+iy$ in $\C^{\,n}$ being of {\em mixed additive type} $(K,v)$, \ie  
\begin{equation*}
  \max\left\{|x_k-\sum_{r=1}^n  j_r x_r|,
\|y_k-\sum_{r=1}^n j_r y_r\|\right\} \geq K|\jb|^{-v}_1 
\end{equation*}
for each $\jb\in (\N\cup \{0\})^n$ with $|\jb|_1\geq 2$ and each 
$k=1,\dots,n$. 

Given $v>0$, the complementary set of points which are not of mixed
additive type for any $K>0$ is null.  It is closely related to and has
the same \HD{} as the simpler set $F_v$ say of points $(x,y)\in\R^{2n}$ such
that
\begin{equation*}
  \max\{|\jb\ip x|,\|\jb\ip y\|\} < |\jb|^{-v} 
\end{equation*}
for infinitely many $\jb\in \Z^n$.  This set is roughly speaking made
up of a `distance from 0' part and a `distance from $\Z$' part.  The
resonant sets are a system of cartesian products of a hyperplane
normal to $\q$ through the origin in $\R^n$ and a family of parallel
hyperplanes.  The lower bound for the dimension is obtained by showing
that the resonant sets form a ubiquitous system, the upper by using
the natural limsup cover for $F_v$~\cite{DRV90a,DRV94}.  If $v \geq
(n-1)/2$, then $ F_v$ and $\cE_v$ are null with Hausdorff dimension
\begin{equation}
\label{eq:dimEv}
\dim F_v=\dim {\cE}_v = 2(n-1) + \frac{n+1}{v+1}.
\end{equation} 

Mixed additive type can also be used in the analysis of linearising
periodic differential equations or, in other words, for finding the
normal form of a vector field on $\C^{\,n}\times \ci$ near a singular
point.  The periodic differential equation
\begin{equation*}
  \dot z = B z+ Q(z,t), 
\end{equation*}
where $B$ is an $n\times n$ complex matrix and $Q\col\C^{\,n}\times \ci\to
\C^{\,n}$ is analytic and has period $2\pi$ in $t$ and satisfies 
\begin{equation*}
Q(0,t)=0, \quad \frac{\partial Q_j(0,t) }{\partial z_k}=0 {\text{ for }}  
1\leq j,k \leq n, 
\end{equation*}
can be linearised to the form $\dot\zeta=B\zeta$ in a neighbourhood of $0$ if
the real and imaginary parts of the eigenvalues of $B$ form a vector
of mixed additive type
(see~\cite[Theorem~B]{DRV89a},~\cite{Siegel52}).  For each $v>0$, the
exceptional set of such vectors has the same \HD{} as
$\cE_v$ in~\eqref{eq:dimEv}.
Similar results hold for autonomous differential equations.

\newcommand{\B}{\mathbb{B}}
\newcommand{\cH}{\cal{H}}
\newcommand{\LG}{\Lambda(G)}
\newcommand{\Lg}{\mu_{g}}
\newcommand{\Mob}{M{\" o}bius}
\newcommand{\Poinc}{Poincar{\'e}} 
\newcommand{\SL}{\text{SL}}
\newcommand{\Sn}{\mathbb{S}^n} 

\newcommand{\x}{\mathbf{x}}

\section[Hyperbolic space]
{\DA{} in hyperbolic space\index{hyperbolic space}}  
\label{sec:dahs}

Diophantine approximation has a fertile interpretation in
hyperbolic space in which many results can be interpreted dynamically,
for example geodesic orbits on manifolds play an important role.  
 The action of a discrete subgroup of the
group of orientation preserving \Mob {} 
transformations~\index{M{\"o}bius transformations} of the upper
half plane to itself allows a fertile generalisation of the classical
theory of \DA.  The rationals $\Q$ can be interpreted as the orbit of
the point at infinity under the linear fractional or \Mob{}
transformation of the extended complex plane $\C\cup\{\infty\}$ in
which
\begin{equation*}
  z\mapsto \frac{az+b}{cz+d}, \ a,b,c,d \in \Z, \ ad-bc=1,
\end{equation*}
\ie as the orbit of $\infty$ under the modular group $\SL(2,\Z)$
acting on points in the upper half plane.  These observations allow a
very nice translation of classical \DA{} into the hyperbolic space
setting by considering the action of a Kleinian group $G$ and the
(Euclidean) distance of the orbit of a special point in the limit set
$\LG$ from the other points in the set. Here $\LG$ is the set of
accumulation points of the orbit $G(x)$ of the point $x$ in hyperbolic
space and turns out to be independent of $x$.  Details of the basic
properties of Kleinian groups acting on hyperbolic space can be found
in~\cite{AhlforsMT,BeardonGDG,NichollsETDG}. From now on, the groups
considered will be non-elementary and geometrically finite.
Approximation of real numbers by rationals is replaced by
approximating points in $\LG$ by points in the orbit of the special point
$\dpt$ which is either a parabolic fixed point of $G$ when $G$ has
them or a hyperbolic fixed point otherwise.  More precisely, when the
unit ball $\B^{n+1}$ endowed with the hyperbolic metric $\rho$ (given
by $d\rho=|dx|_2/(1-|x|_2^2)$) is the model for hyperbolic space,
given a point $\xi\in\LG$, we consider the quantity
 \begin{equation}
|\xi-g({\dpt})|_2 < \lambda_g^{-v}
   \label{eq:hypDA}
 \end{equation}
 as $g$ runs through $G$.  
The 
 `denominator' in this  setting is the reciprocal of the
 modulus of the Jacobian $g'(0)$ of each $g\in G$ at $0$,
\begin{equation*}
\lag = |g'(0)|^{-1} 
\end{equation*}
and is comparable to $e^{\rho(0,g(0))}$ (\ie $\lag e^{-\rho(0,g(0))}$ is
bounded above and below by positive  constants). In the upper half plane
model $\lag=|g'(i)|^{-1}$.  Thus $\lag\to\infty$ as $|g(0)|\to 1$,
\ie as the orbit of the origin moves out towards the boundary of
hyperbolic space~\cite{NichollsETDG}.

It turns out that not only is $\lag$  an 
appropriate analogue for the modulus of
the denominator of a rational $p/q$ in the classical setting 
but also properties of the Dirichlet series $\sum_{g\in G} \lag^{-s}$
are  connected with the \HD{} of $\LG$.  In fact 
the {\em exponent of convergence}
\begin{equation*}
\d (G) = \inf \{s > 0\col \sum_{g \in G} \lag^{-s} < \infty \} 
 \end{equation*}
 of $G$ is equal to $\dim \LG$~\cite{SV95,Sullivan84}.  The
 appropriateness of $\lag$ as a `denominator' of $g$ can be seen
 clearly by considering the case $G=\SL(2,\Z)$ and $\dpt=\infty$.
 Each group element $g$ satisfies $g(\infty)=p/q$ for some $p/q\in\Q$
 and a straight\-forward calculation yields that $\lag$ is comparable
 to $q^2$. By
 analogy with the classical case, a point $\xi\in \LG$ is called {\em
   $v$-approximable} if the inequality
\begin{equation*}
|\xi-g(\dpt)|_2<{\lambda_g}^{-v}
\end{equation*}
holds for infinitely many $g\in G$.

When $G$ has parabolic fixed points (the most interesting case),
Dirichlet's theorem has the following analogue: There exists a
constant $C>0$ depending only on $G$, such that for each $\xi\in \LG$
and each integer $N>1$,
\begin{equation*}
  |\xi-g(\dpt)|_2< C/\sqrt{\lag N }\le C/\lag 
\end{equation*}
for some parabolic $\dpt$ and $g\in G$ with $\lag\leq
N$~\cite{Patterson76a,VelaniQF}.  The analogy can be pursued further,
in the first place to Khintchine's theorem~\cite{Patterson76a,SV95}. We suppose $G$
has a parabolic element $\dpt$ of rank $r_{\dpt}$.  Let $\mu$ be
Patterson measure~\cite{Patterson76b}~\index{Patterson measure} and
$\psi\col [1/2,\infty) \to \R^+$ be decreasing and satisfy
$\psi(2x)>c\psi(x)$ for some constant $c>0$.  Then the set
$(W(G,\dpt;\psi))$ of points $\xi$ satisfying
\begin{equation*}
  |\xi-g(\dpt)|_2< \psi(\lag) 
\end{equation*}
infinitely often has Patterson measure 0 or $1$ according as the sum
\begin{equation*}
  \sum_{k=1}^\infty \psi(K^k)^{2\delta(G)-r_{\dpt}}
\end{equation*}
converges or diverges for some $K>1$ (when $G$ has no parabolic
elements, the sum is simpler).

The set $W_v(\SL(2,\Z),\infty)$ of $v$-approximable points in the limit set  
corresponds to $\K_v$, the set of $v$-approximable points in $\R$,
although the exponent $v$ is normalised differently.  For any $\e>0$,
\begin{equation*}
W_v(\SL(2,\Z),\infty) \subset \K_{2v+1} \subset W_{v-\e}(\SL(2,\Z),\infty),
\end{equation*} 
so that $\dim W_v(\SL(2,\Z),\infty) = 1/(v+1)$ by the \JB{} theorem.
Moreover a more general counterpart of the 
theorem holds~\cite{HVJBGF}. 
\begin{theorem}
  Let $G$ be a non-elementary geometrically finite group and suppose
  $v\geq 0$.  If $G$ has a parabolic element $\dpt$ say, then
  \begin{equation*}
    \dim W_v(G,\dpt)=\min
\left\{\frac{\delta(G)+r_{\dpt}}{2v+1},\frac{\delta(G)}{v+1}\right\},
  \end{equation*}
where $r_{\dpt}$ is the rank of $\dpt$. 
Otherwise $G$ is convex co-compact and 
$$
\dim  W_v(G,\dpt)= \delta(G)/(v+1).
$$
\end{theorem}

This setting provides a beautiful dynamical interpretation of
approximation.  When $G$ is \gf, the quotient space $M=\B^{n+1}/G$
obtained by identifying equivalent points in $\B^{n+1}$ under the
action of $G$ is an $n$-\dm al manifold.  When $G$ is a Fuchsian
group, the manifold is a Riemann surface of constant negative
curvature, when $n=1$ and $G=\SL(2,\Z)$, $\Hy^2/\SL(2,\Z)$, where
$\Hy^2$ is the upper half plane, is the modular surface.  The manifold
$M$ can be decomposed into a disjoint union of a compact part with a
finite number of exponentially `narrowing' cuspidal ends
(corresponding to a set of inequivalent parabolic fixed points of $G$)
and `exploding' ends or funnels (corresponding to the free faces of a
convex fundamental polyhedron for $G$).  The set of points $\xi$ in
$\LG$ for which there exists a constant $c(\xi)$ such that
\begin{equation*}
  |\xi-g(\dpt)|\geq c(\xi)/\lag
\end{equation*} 
holds for all $g\in G$ is the set of badly approximable points,
denoted by $\BA(G,\dpt)$. When $\dpt$ is a parabolic fixed point of
$G$, the set $\BA(G,\dpt)$ corresponds to `bounded' geodesics on the
manifold; the dimension results for $\BA(G,\dpt)$ imply that $M$ is of
full Hausdorff dimension $\d(G)$~\cite{BJ97,FM95}.  On the other hand,
the set $W(G,\dpt;\psi)$ corresponds to excursions by geodesics into
the cuspidal ends of the manifold (the rate being governed by the
function $\psi$); \ie to `divergent' geodesics which persistently
enter a `shrinking' neighbourhood of the cusp associated with $\dpt$.
The work of Dani and and Margulis on the structure of bounded and
divergent trajectories in homogeneous spaces has had a profound
influence on Diophantine
approximation~\cite{Dani84,Dani85,Dani86a,Dani86c,Dani86b,DaniNTDS,DM93,Margulis87}.

 The  notion of a shrinking target\index{shrinking target} introduced
in~\cite{HV95} arose from a general consideration of Diophantine
approximation and has proved fruitful; it is exploited in~\cite{HVSTT} 
and~\cite{HVMDA}.  Let $X$ be a metric space with a Borel
probability measure $\mu$ and $T\col X\to X$ measure preserving and 
ergodic\index{ergodic}.  Then given the function $\varrho\col
\R^+\to\R^+$, the set
\begin{equation*}
 W(a, \varrho)=  \{x\in X\col T^q(x)\in B(a, \varrho(q)) 
{\text{ for infinitely many }} q\in \N\}, 
\end{equation*}
where $B(a,\e)$ is a ball of radius $\e$ in $X$, is an analogue of the
sets $W(G,\dpt;\psi)$ and $\K^{(1,n)}(\psi)$.  By the Birkhoff ergodic
theorem~\cite{KatokHass}\index{Birkhoff ergodic theorem}, if
$\varrho(q)=\varrho_0$, a constant, for $q$ sufficiently large and
$\mu(B(a, \varrho))>0$, the set
$$
  \{x\in X\col T^q(x)\in B(a,\varrho_0) 
{\text{ for infinitely many }} q\in \N\}
$$
has full $\mu$-measure, in other words the forward orbit of almost
all points falls into the ball infinitely often.  It follows that the
complement of $W(a, \varrho_0)$ in $X$, consisting of points whose
forward orbits land in the ball only a finite number of times, is of
zero $\mu$-measure. This set corresponds to the set of badly
approximable points in a dynamical system which has been studied 
 in~\cite{DaniNTDS,Urb91}.

On the other hand,  when
$ \varrho(q)\to 0$ as $q\to\infty$, points in $W(a, \varrho)$
have trajectories which hit a shrinking ball or target infinitely
often and  are called {\em $\varrho$-approximable}.  
 Points in the backward orbit of $a$ in $X$ are resonant
points corresponding  to the
rationals in  $\K^{(1,n)}_v$ and orbit points in $W(G,\dpt;\psi)$.  
In this framework,  the theory of \DA{} in the hyperbolic
space setting can be used  with ideas from ergodic theory to
analyse the structure of a variety of
apparently quite unrelated sets in complex dynamics and other 
dynamical systems.  For further details see~\cite{HVMDA}.

Expanding Markov maps $T$ of the unit interval (such as the continued
fraction map) have been analysed along these lines~\cite{an97b}
and the theory has been extended to higher dimensional tori and to
maps which are multiplication by integer matrices modulo
$\Z^n$~\cite{HVSTT}.  The very well approximable sets associated with
a given dynamical system have unexpected links with exceptional sets
arising from points in the phase space which have `badly behaved'
ergodic averages and with multifractal spectra~\cite{Falc98,HVMDA}.

\section{Extremal manifolds and flows}
\label{sec:flows}

So far, we have been looking at sets which arise in dynamical systems
through Diophantine conditions and using number theory to show that
these sets are exceptional (\ie null) and determining their \HD.
However, the tables have been turned with dynamical systems ideas
being used to great effect, first to prove Oppenheim's
conjecture~\cite{Margulis87} and more recently in Diophantine
approximation, including the important topic of approximation on
manifolds~\cite{Kleinbock01}.  The latter stemmed from Mahler's
conjecture~\index{Mahler's conjecture} work in transcendence theory
during the 1930's.  He conjectured that the set of points
$(x,x^2,\dots,x^n)$ on
the Veronese curve which are  $v$-simultaneously approximable 
is relatively
null for $v>1/n$, \ie
\begin{equation*}
  |\{x\in\R\col \max_{ j=1,\dots,n}\left\{\|qx^j\|\right\} 
   <q^{-v} {\text{ for infinitely many }} q\in \N\}|=0
\end{equation*}
when $v>1/n$.  This property, which is a special case of \Kh's theorem
and says that the exponent in Dirichlet's theorem on simultaneous
Diophantine approximation cannot be improved for almost all points, is
called {\em extremality}\index{extremality}.  By \Kh's transference
principle, the set of simultaneously VWA points is also dually VWA, in
the sense that the set of points $\xi$ satisfying
\begin{equation*}
 |q_0+q_1\xi+\dots+q_n\xi^n|<|\q|^{-v} 
\end{equation*}
for infinitely many $\q \in\Z^n$ is  null for $v>n$.
Mahler's conjecture was proved by
\Spr~\cite{Spr66} who then conjectured that non-degenerate analytic
manifolds were extremal.  Extremality was later strengthened to a
stronger `multiplicative' form by Baker~\cite{ABakerTNT} and
\Spr~\cite{Spr80} in which given any $v>1$, the inequality 
\begin{equation*}
  \|\q\ip y\|\leq \prod_{j=1}^n\max\{1,|q_j|\}^{-v}
\end{equation*}
holds for infinitely many $q\in \N$ for almost all points
$y=(y_1,\dots,y_n)\in M$. Such points are called {\em very well
  multiplicatively approximable} or VWMA\index{VWMA}.  The `strong
extremality'\index{strong extremality} conjecture was proved by
Kleinbock and Margulis~\cite{km98a} for smooth manifolds which are
non-degenerate (or `non-flat' locally) almost everywhere.  Their proof
uses actions on the lattice
\begin{equation*}
  \Lambda_{y}=\left(\begin{matrix} 1 & y^T\\
                            0 & I_n \end{matrix}  \right) \, \Z^{n+1},
\end{equation*}
where $y\in\R^n$, in the homogeneous space
$\SL_{n+1}(\R)/\SL_{n+1}(\Z)$ of unimodular lattices by the semigroup 
$\{g_{\tv} \col \tv\in \N^n\}$, where
\begin{equation*}
  g_{\tv}=\diag (e^{\sum_{j=1}^n t_j}, e^{-t_1},\dots,e^{-t_n}).
\end{equation*}
Elements in $ \Lambda_{y}$ are of the form $(q_0+\q\ip y,\q)$, $\q\in
\Z^n$, $q_0\in\Z$ and the distance of the lattice $\Lambda$ from the
origin is
\begin{equation*}
  \d(\Lambda)=\inf_{a\in\Lambda\setminus 0} |a|
\end{equation*} 
(recall $|a|$ is the height of $a$).   By Mahler's Compactness
Criterion, an infinite sequence $\Lambda^{(k)}$, $k=1,2,\dots,$
diverges to infinity in $\SL_{n+1}(\R)/\SL_{n+1}(\Z)$ if and only if $
\d(\Lambda^{(k)})\to 0$ as $k\to\infty$. It can be shown that if a
point $y$ is very well multiplicatively approximable, then there
exists a $\gamma>0$ such that $\d(g_{\tv}\Lambda_{y})\leq e^{-\gamma
  |\tv|_1}$ for infinitely many $\tv\in \N^n$.  By the Borel-Cantelli
Lemma, the set of VWMA points on the manifold is exceptional (\ie the
manifold is strongly extremal) if the sum
$$
\sum_{\tv\in\N^n}|\{\x\in B\col \d(g_{\tv}\Lambda_{f(x)})\leq e^{-\gamma
  |\tv|_1}|<\infty,
$$
where $B$ is a neighbourhood of a non-degenerate point $x_0$ and
$y=f(x)$, where $f$ is the local parametrisation function.  This is
established by modifying the arguments of Dani and
Margulis~\cite{Dani86c}.

\section{Conclusion}

Exceptional sets in dynamical systems can arise from Diophantine
conditions and can be analysed using number theory. The number theory
has in turn served as a model for an analogous theory in groups
actions in hyperbolic space.  The interpretation of the orbit
structure in terms of Diophantine approximation has been very
fruitful, leading to a general notion of shrinking targets and in
another dramatic development to the recent solution of the
Baker-\Spr{}\index{Baker-Sprind{\v z}uk conjecture} conjecture and the
proof of multiplicative Khintchine-type theorems for $C^n$ manifolds
in $\R^n$~\cite{bkm01}.  However, Beresnevich's proof~\cite{vvb01},
which uses ideas based on those of \Spr, of a \Kh-type theorem in the
case of convergence for manifolds requires weaker differentiability
conditions, showing that classical methods remain effective.  The
complementary case of divergence is much more difficult but progress
is being made.

\section{Acknowledgments}

Thanks are due to the organisers of the Ergodic Theory, Geometric
Rigidity and Number Theory Conference and the Newton Institute for
arranging an interesting and varied programme and for their
hospitality and to Marc Burger and Alessandra Iozzi for organising this
volume and for their patience.  I am also grateful to Detta Dickinson
for very helpful discussions and to the referee whose corrections and
suggestions have improved this paper. 

\bibliographystyle{amsplain}

\providecommand{\bysame}{\leavevmode\hbox to3em{\hrulefill}\thinspace}


\begin{thebibliography}{10}

\bibitem{an97b}
A.~G. Abercrombie and R.~Nair, \emph{An exceptional set in the ergodic theory
  of {M}arkov maps on the interval}, Proc.\ Lond.\ Math.\ Soc. \textbf{75}
  (1997), 221--240.

\bibitem{AhlforsMT}
L.~V. Ahlfors, \emph{M{\"o}bius transformations in several dimensions}, Lecture
  Notes, School of Mathematics, University of Minnesota, 1954.

\bibitem{Arnold63}
V.~I. Arnol'd, \emph{Small denominators and problems of stability of motion in
  classical and celestial mechanics}, Usp.\ Mat.\ Nauk \textbf{18} (1963),
  91--192, English transl. in Russian Math.\ Surveys, {\bf{18}}\, (1963),
  85--191.

\bibitem{ArnoldMM}
\bysame, \emph{Mathematical methods of classical mechanics}, Springer-Verlag,
  1978, Translated by K.~Vogtmann and A.~Weinstein.

\bibitem{ArnoldGM}
\bysame, \emph{Geometrical methods in ordinary differential equations},
  Springer-Verlag, 1983, Translated by J.~Sz{\"u}cs.

\bibitem{AKNIII}
V.~I. Arnol'd, V.~V. Kozlov, and A.~I. Neishtadt, \emph{Mathematical aspects of
  classical and celestial mechanics}, Encyclopaedia of Mathematical Sciences,
  vol.\ 3, Dynamical Systems III, Springer-Verlag, 1980, Translated by
  A.~Jacob.

\bibitem{APIDS}
D.~K. Arrowsmith and C.~M. Place, \emph{An introduction to dynamical systems},
  Cambridge University Press, 1990.

\bibitem{ABakerTNT}
A.~Baker, \emph{Transcendental number theory}, second ed., Cambridge University
  Press, 1979.

\bibitem{BS}
A.~Baker and W.~M. Schmidt, \emph{Diophantine approximation and {H}ausdorff
  dimension}, Proc.\ Lond.\ Math.\ Soc. \textbf{21} (1970), 1--11.

\bibitem{BeardonGDG}
A.~F. Beardon, \emph{The geometry of discrete groups}, Springer-Verlag, 1983.

\bibitem{vvb01}
V.~V. Beresnevich, \emph{A {G}roshev type theorem for convergence on
  manifolds}, Acta Math.\ Hung. (2001), to appear.

\bibitem{MDAM}
V.~I. Bernik and M.~M. Dodson, \emph{Metric {D}iophantine approximation on
  manifolds}, Cambridge Tracts in Mathematics, No.~137, Cambridge University
  Press, 1999.

\bibitem{bkm01}
V.~I. Bernik, D.~Y. Kleinbock, and G.~A. Margulis, \emph{Khintchine-type
  theorems for manifolds: the convergence case for standard and multiplicative
  versions}, submitted to Inter.\ Math.\ Res.\ Notices.

\bibitem{Bes34}
A.~S. Besicovitch, \emph{Sets of fractional dimensions ({I}{V}): on rational
  approximation to real numbers}, J.\ Lond.\ Math.\ Soc. \textbf{9} (1934),
  126--131.

\bibitem{BJ97}
C.~J. Bishop and P.~W. Jones, \emph{Hausdorff dimension and {K}leinian groups},
  Acta Math. \textbf{111} (1997), 1--39.

\bibitem{BD86}
J.~D. Bovey and M.~M. Dodson, \emph{The {H}ausdorff dimension of systems of
  linear forms}, Acta Arith. \textbf{45} (1986), 337--358.

\bibitem{Cassels}
J.~W.~S. Cassels, \emph{An introduction to {D}iophantine approximation},
  Cambridge Tracts in Math.\ and Math.\ Phys., No.~45, Cambridge University
  Press, 1957.

\bibitem{Dani84}
S.~G. Dani, \emph{On orbits of unipotent flows on homogeneous spaces}, Ergod.\
  Th.\ Dyn.\ Sys. \textbf{4} (1984), 25--34.

\bibitem{Dani85}
\bysame, \emph{Divergent trajectories of flows on homogeneous spaces and
  homogeneous {D}iophantine approximation}, J.\ reine angew.\ Math.
  \textbf{359} (1985), 55--89.

\bibitem{Dani86a}
\bysame, \emph{Bounded orbits of flows on homogeneous spaces}, Comm.\ Math.\
  Helv. \textbf{61} (1986), 636--660.

\bibitem{Dani86c}
\bysame, \emph{On orbits of unipotent flows on homogeneous spaces, {II}},
  Ergod.\ Th.\ Dyn.\ Sys. \textbf{6} (1986), 167--182.

\bibitem{Dani86b}
\bysame, \emph{Orbits of horospherical flows}, Duke Math.\ J. \textbf{53}
  (1986), 177--188.

\bibitem{DaniNTDS}
\bysame, \emph{On badly approximable numbers, {S}chmidt games and bounded
  orbits of flows}, Number theory and dynamical systems (M.~M. Dodson and
  J.~A.~G. Vickers, eds.), LMS Lecture Note Series, vol. 134, Cambridge
  University Press, 1987, pp.~69--86.

\bibitem{DM93}
S.~G. Dani and G.~A. Margulis, \emph{Limit distributions of orbits of unipotent
  flows and values of quadratic forms}, Adv.\ Soviet Math. 16, Amer.\ Math.\
  Soc., RI \textbf{16} (1993), 91--137.

\bibitem{KAMtut}
R.~de~la Llave, \emph{A tutorial on {KAM} theory}, Smooth Ergodic Theory and
  Applications, Proceedings of Symposia in Pure Mathematics, Summer Research
  Institute, Seattle WA, Amer.\ Math.\ Soc., 1999, to appear 2001.

\bibitem{HD93}
H.~Dickinson, \emph{The {H}ausdorff dimension of systems of simultaneously
  small linear forms}, Mathematika \textbf{40} (1993), 367--374.

\bibitem{dgy95}
H.~Dickinson, T.~Gramchev, and M.~Yoshino, \emph{First order pseudodifferential
  operators on the torus: normal forms, {D}iophantine approximation and global
  hypoellipticity}, Ann.\ Univ.\ Ferrara, Sez.\ VII -- Sc.\ Mat. \textbf{61}
  (1995), 51--64.

\bibitem{mmd92}
M.~M. Dodson, \emph{Hausdorff dimension, lower order and {K}hintchine's theorem
  in metric {D}iophantine approximation}, J.\ reine angew. Math. \textbf{432}
  (1992), 69--76.

\bibitem{DRV89a}
M.~M. Dodson, B.~P. Rynne, and J.~A.~G. Vickers, \emph{Averaging in
  multi-frequency systems}, Nonlinearity \textbf{2} (1989), 137--148.

\bibitem{DRV90a}
\bysame, \emph{Diophantine approximation and a lower bound for {H}ausdorff
  dimension}, Mathematika \textbf{37} (1990), 59--73.

\bibitem{DRV94}
\bysame, \emph{The {H}ausdorff dimension of exceptional sets associated with
  normal forms}, J.\ Lond.\ Math.\ Soc. \textbf{49} (1994), 614--624.

\bibitem{DV}
M.~M. Dodson and J.~A.~G. Vickers, \emph{Exceptional sets in
  {K}olmogorov-{A}rnol'd-{M}oser theory}, J.\ Phys.\ A \textbf{19} (1986),
  349--374.

\bibitem{NTDS}
M.~M. Dodson and J.~A.~G. Vickers~(eds.), \emph{Number theory and dynamical
  systems}, Lond.\ Mathematical Society Lecture Note Series, vol. 134,
  Cambridge University Press, 1989.

\bibitem{FalcGFS}
K.~Falconer, \emph{The geometry of fractal sets}, Cambridge Tracts in
  Mathematics, No.~85, Cambridge University Press, 1985.

\bibitem{Falc98}
\bysame, \emph{Representation of families of sets, dimension spectra and
  {D}iophantine approximation}, Math.\ Proc.\ Cam. Philos.\ Soc. \textbf{128}
  (2000), 111--121.

\bibitem{FM95}
J.-L. Fernandez and M.~V. Meli{\'a}n, \emph{Bounded geodesics of {R}iemann
  surfaces and hyperbolic manifolds}, Trans.\ Amer.\ Math.\ Soc. \textbf{9}
  (1995), 3533--3549.

\bibitem{HW}
G.~H. Hardy and E.~M. Wright, \emph{An introduction to the theory of numbers},
  4th ed., Clarendon Press, 1960.

\bibitem{HarmanMNT}
G.~Harman, \emph{Metric number theory}, LMS Monographs New Series, vol.~18,
  Clarendon Press, 1998.

\bibitem{Herman87}
M.~R. Herman, \emph{Recent results and some open questions on {S}iegel's
  linearisation theorem on germs of complex analytic diffeomorphisms of
  {$\C^n$} near a fixed point}, Proceedings of the VIIIth International
  Congress of Mathematical Physics, Marseilles 1986 (M.~Mebkhout and
  R.~S{\/e}n{\/e}or, eds.), World Scientific, 1987, pp.~138--184.

\bibitem{HV95}
R.~Hill and S.~L. Velani, \emph{Ergodic theory of shrinking targets}, Invent.\
  Math. \textbf{119} (1995), 175--198.

\bibitem{HVMDA}
\bysame, \emph{Metric {D}iophantine approximation in {J}ulia sets of expanding
  rational maps}, Publ.\ Math.\ I.H.E.S. (1997), no.~85, 193--216.

\bibitem{HVJBGF}
\bysame, \emph{The {J}arn{\'\i}k-{B}esicovitch theorem for geometrically finite
  {K}leinian groups}, Proc.\ Lond.\ Math.\ Soc. \textbf{77} (1998), 524--550.

\bibitem{HVSTT}
\bysame, \emph{The shrinking target problem for matrix transformations of
  tori}, J.\ Lond.\ Math.\ Soc. \textbf{60} (1999), 381--398.

\bibitem{Ja28}
V.~Jarn{\'\i}k, \emph{Zur metrischen {T}heorie der diophantischen
  {A}pproximationen}, Prace Mat.-Fiz. (1928-9), 91--106.

\bibitem{Ja29}
\bysame, \emph{Diophantischen {A}pproximationen und {H}ausdorff\-sches {M}ass},
  Mat.\ Sbornik \textbf{36} (1929), 371--382.

\bibitem{Ja31}
\bysame, \emph{{\"U}ber die simultanen diophantischen {A}pproximat\-ionen},
  Math.\ Z. \textbf{33} (1931), 503--543.

\bibitem{KatokHass}
A~Katok and B.~Hasselblatt, \emph{Introduction to the modern theory of
  dynamical systems}, Cambridge University Press, 1995.

\bibitem{Kleinbock01}
D.~Y. Kleinbock, \emph{Applications of ergodic theory to metric {D}iophantine
  approximation}, Smooth Ergodic Theory and Applications, Proceedings of
  Symposia in Pure Mathematics, Summer Research Institute, Seattle WA, American
  Mathematical Society, 1999, to appear 2001.

\bibitem{Kleinbock99}
\bysame, \emph{Badly approximable systems of affine forms}, J.~Number Th.
  \textbf{79} (1999), 83--102.

\bibitem{km98a}
D.~Y. Kleinbock and G.~A. Margulis, \emph{Flows on homogeneous spaces and
  {D}iophantine approximation on manifolds}, Ann.\ Math. \textbf{148} (1998),
  339--360.

\bibitem{Kolmogorov54}
A.~N. Kolmogorov, \emph{On the conservation of conditionally periodic motions
  under small perturbations of the {H}amiltonian}, Dokl.\ Akad.\ Nauk SSSR
  \textbf{98} (1954), 527--530, in Russian.

\bibitem{JL98}
J.~Levesley, \emph{A general inhomogeneous {\JB}{} theorem}, J.~Number Th.
  \textbf{71} (1998), 65--80.

\bibitem{Margulis87}
G.~A. Margulis, \emph{Formes quadratiques ind{\'e}finies et flots unipotents
  sur l'espaces homog{\`e}nes}, C.\ R.\ Acad.\ Sci.,\ Paris, Ser.\ I
  \textbf{304} (1987), 249--253.

\bibitem{MattilaGS}
P.~Mattila, \emph{Geometry of sets and measures in {E}uclidean space},
  Cambridge University Press, 1995.

\bibitem{Moser62}
J.~Moser, \emph{On invariant curves of area-preserving maps of an annulus},
  Nachr.\ Akad.\ Wiss.\ G{\"o}tt., Math.\ Phys.\ Kl. (1962), 1--20.

\bibitem{MoserSRM}
\bysame, \emph{Stable and random motions in dynamical systems}, Princeton
  University Press, 1973.

\bibitem{Moser78}
\bysame, \emph{Is the solar system stable?}, Math.\ Intelligencer \textbf{1}
  (1978), 65--71.

\bibitem{NichollsETDG}
P.~J. Nicholls, \emph{The ergodic theory of discrete groups}, LMS Lecture
  Notes, vol. 143, Cambridge University Press, 1989.

\bibitem{NiteckiDD}
Z.~Nitecki, \emph{Differentiable dynamics}, MIT\ Press, 1971.

\bibitem{Patterson76a}
S.~J. Patterson, \emph{Diophantine approximation in {F}uchsian groups}, Phil.\
  Trans.\ Roy.\ Soc.\ Lond.\ A \textbf{262} (1976), 527--563.

\bibitem{Patterson76b}
\bysame, \emph{The limit set of a {F}uchsian group}, Acta Math. \textbf{136}
  (1976), 241--273.

\bibitem{Pesin97}
Ya.~B. Pesin, \emph{Dimension theory in dynamical systems. {C}ontemporary views
  and applications}, Chicago Lectures in Math., University of Chicago
  Press, 1997.

\bibitem{Poschel82}
J.~P{\"o}schel, \emph{Integrability of {H}amiltonian systems on {C}antor sets},
  Comm.\ Pure Appl. \ Math. \textbf{35} (1982), 653--696.

\bibitem{Rogers}
C.~A. Rogers, \emph{Hausdorff measure}, Cambridge University Press, 1970.

\bibitem{Russmann83}
H.~R. R{\"u}ssmann, \emph{On the existence of invariant curves of twist
  mappings of the annulus}, Geometric {D}ynamics, Lecture Notes in Mathematics,
  vol. 1007, Springer-Verlag, 1983, pp.~677--712.

\bibitem{Russmann91}
\bysame, \emph{On the frequencies of quasi-periodic solutions of analytic
  nearly integrable {H}amiltonian systems}, Progress in {N}onlinear
  {D}ifferential {E}quations and their applications, 12 (V.~Lazutkin S.~Kuksin
  and J.~P{\"o}schel, eds.), vol.~12, Birkh{\"a}user Verlag, Basel, 1994,
  pp.~51--58.

\bibitem{Schmidt69}
W.~M. Schmidt, \emph{Badly approximable systems of linear forms}, J.\ Number
  Th. \textbf{1} (1969), 139--154.

\bibitem{Siegel42}
C.~L. Siegel, \emph{Iteration of analytic functions}, Ann.\ Math.\textbf{43}
  (1942), 607--612.

\bibitem{Siegel52}
\bysame, \emph{{\"U}ber die {N}ormalform analytischer {D}ifferentialgleichungen
  in der {N}{\"a}he einer {G}leichgewichtsl{\"o}sung}, Nachr.\ Akad.\ Wiss.\
  G{\"o}tt.\ Math-Phys.\ Kl (1952), 21--30.

\bibitem{SiegelMoserCM}
C.~L. Siegel and J.~K. Moser, \emph{Lectures on celestial mechanics},
  Springer-Verlag, 1971.

\bibitem{Spr66}
V.~G. Sprind{\v z}uk, \emph{A proof of {M}ahler's conjecture on the measure of
  the set of {$S$}-numbers}, Amer.\ Math. \ Soc.\ Transl.\ Ser. 2 \textbf{51}
  (1966), 215--272.

\bibitem{Sprindzuk}
\bysame, \emph{Metric theory of {D}iophantine approximations}, John Wiley,
  1979, Translated by R.~A.~Silverman.

\bibitem{Spr80}
\bysame, \emph{Achievements and problems in {D}iophantine approximation
  theory}, Usp.\ Mat.\ Nauk \textbf{35} (1980), 3--68, English transl. in
  Russian Math.\ Surveys, {\bf{35}}\, (1980), 1--80.

\bibitem{SV95}
B.~Stratmann and S.~L. Velani, \emph{The {P}atterson measure for geometrically
  finite groups with parabolic elements, new and old}, Proc.\ Lond.\ Math. \
  Soc. \textbf{71} (1995), 197--220.

\bibitem{Sullivan84}
D.~Sullivan, \emph{Entropy, {H}ausdorff measures old and new, and the limit set
  of geometrically finite {K}leinian groups}, Acta Math. \textbf{153} (1984),
  259--277.

\bibitem{Urb91}
M.~Urb{\'a}nski, \emph{The {H}ausdorff dimension of the set of points with
  non-dense orbit under a hyperbolic dynamical system}, Nonlinearity \textbf{4}
  (1991), 385--397.

\bibitem{VelaniQF}
S.~L. Velani, \emph{An application of metric {D}iophantine approximation in
  hyperbolic space to quadratic forms}, Publ.\ Math. \textbf{38} (1994),
  175--185.

\bibitem{Yoccoz92}
J.~C. Yoccoz, \emph{An introduction to small divisors problems}, From number
  theory to physics, Springer-Verlag, 1992, Les Houches, 1989, pp.~659--679.

\bibitem{Yoccoz95a}
\bysame, \emph{Petits diviseurs en dimension 1}, Ast{\'e}risque \textbf{231}
  (1995).

\end{thebibliography}

\printindex

\end{document}